\documentclass[11pt]{article}
\usepackage{amssymb}
\newtheorem{lem}{Lemma}[section]

\newtheorem{pro}{Proposition}[section]
\newtheorem{thm}{Theorem}[section]

\makeatletter\@addtoreset{equation}{section}
\renewcommand\theequation{\thesection.\@arabic\c@equation}

\newenvironment{Literature}[1]
{\begin{thebibliography}{#1}}{\end{thebibliography}}

\begin{document}

\begin{center}
{\LARGE \bf  Dirac eigenvalues  estimates   }   \\

\bigskip

{\LARGE \bf  in terms of symmetric tensors }

\bigskip  \bigskip
{\large  Eui Chul Kim}
\end{center}

\bigskip \bigskip \noindent
{\bf Abstract:} We review some recent results concerning lower
eigenvalues estimates for the Dirac operator [6, 7].
 We show that Friedrich's inequality can be improved via certain well-chosen
symmetric tensors and provide an application
 to Sasakian spin manifolds.

\noindent
\section{Introduction }

\noindent Let $(M^n, g), \,  n \geq 3 , $ be an $n$-dimensional closed
Riemannian spin manifold. Let $(E_1, \ldots, E_n )$ be  a local
orthonormal frame field on $(M^n, g)$. Then the spinor derivative
$\nabla$ and the Dirac operator $D$, acting on sections $\psi \in
\Gamma( \Sigma(M^n))$ of the spinor bundle $\Sigma(M^n)$ over $(M^n,
g)$, are locally expressed as
\begin{eqnarray*}
\nabla_X \psi = X(\psi) + \frac{1}{4} \sum_{i=1}^n E_i \cdot
\nabla_X E_i \cdot \psi \ , \quad D \psi = \sum_{i=1}^n E_i \cdot
\nabla_{E_i} \psi  ,
\end{eqnarray*}
 respectively, where the dot "$\cdot $" indicates the Clifford multiplication [2, 4].
 Friedrich proved in [3] that the
smallest eigenvalue $\lambda_1$ of $D$  satisfies
\begin{equation}
\lambda_1^2    \, \geq  \, \frac{n}{4(n-1)} \, {\rm inf}_M \, S   ,
\end{equation}
where $S$ is the scalar curvature of $(M^n, g)$. The limiting case
of (1.1) occurs if and only if $(M^n, g)$ admits a nontrivial spinor
field $\psi$ called {\it Killing spinor}, satisfying
\begin{equation}
\nabla_X \psi =  - \frac{\lambda_1}{n}  \, X  \cdot \psi ,
\end{equation}
where $X$ is an arbitrary vector field on $M^n$. The
simply-connected manifolds $(M^n, g)$  admitting Killing spinors
were classified by B\"{a}r [1], namely, the limiting manifold $(M^n,
g)$ must be either a standard n-sphere, an Einstein-Sasaki
 manifold, a  6-dimensional nearly K\"{a}hler
  manifold or a 7-dimensional
 manifold with 3-form $\phi$, $\nabla \phi = \ast \phi$. Note that all of these limiting manifolds are Einstein, since equation (1.2) allows a nontrivial solution $\psi$ only if $(M^n,
g)$ is Einstein.

\bigskip
It has been found that inequality (1.1) is not optimal if $(M^n , g)$ allows certain geometric structures,
 since the limiting case of (1.1) can not be
attained [9, 10, 11]. For example, Kirchberg
proved for K\"{a}hler spin manifolds that  the
smallest eigenvalue $\lambda_1$ of the Dirac operator satisfies
\begin{equation}    \lambda_1^2  \  \geq \  \frac{n+2}{4n} \, {\rm inf}_M \,
S   \qquad   \mbox{for}  \quad  n \equiv 2 \, \mbox{mod} \, 4
\end{equation}
and
\begin{equation}    \lambda_1^2  \  \geq \  \frac{n}{4(n-2)} \, {\rm inf}_M \,
S    \qquad    \mbox{for}  \quad n \equiv 0 \, \mbox{mod} \, 4 .
\end{equation}

\bigskip
Improvements of Friedrich's inequality (1.1) do typically depend on
additional geometric structures on the considered manifold $(M^n,
g)$. The aim of this article is to review some new results in [6, 7], showing that
Friedrich's inequality can be improved via divergencefree symmetric tensors  as well as Codazzi tensors (see Theorem 2.1 and 2.2).
In the last section we discuss geometric implications of Theorem 2.1  over Sasakian spin manifolds.

\noindent
\section{Dirac eigenvalues estimates  in terms of symmetric tensors }

\noindent
Throughout the article we fix some terminology.

\bigskip \noindent
{\bf Definition 2.1} Let $P$ be a first order self-adjoint elliptic
operator on some closed Riemannian spin manifold. An eigenvalue
$\lambda \in {\mathbb R}$ of $P$ is called the {\it first
eigenvalue}  if $\lambda^2$ is the smallest eigenvalue of $P^2$. An
eigenspinor $\varphi$ of $P$ is called a {\it  first eigenspinor} if
its associated eigenvalue $\lambda$ is the first eigenvalue of $P$.

\bigskip
\noindent Evidently, the Killing spinors satisfying equation (1.2) are first eigenspinors
of the Dirac operator. Let's see one more example.
The limiting case of inequality (1.3) occurs if and only if the coupled system
\begin{eqnarray}
D \psi & = & \lambda_1 \psi  ,    \nonumber  \\
&      &      \nonumber      \\
 \nabla_X  \psi & = & - \frac{\lambda_1}{n+2} \, X \cdot  \psi
+ \frac{\lambda_1}{n+2} \, J(X) \cdot \Omega \cdot \psi
\end{eqnarray}
admits a nontrivial solution $\psi$ called {\it K\"{a}hlerian Killing spinor}, where
$\Omega$ is the K\"{a}hler form.  Thus, the  K\"{a}hlerian Killing spinors are first eigenspinors
of the Dirac operator.

\bigskip
Let us now consider a nondegenerate symmetric (0,
2)-tensor field $\beta$ on $(M^n, g)$ and  define the $\beta$-{\it
twist} $D_{\beta}$ of the Dirac operator $D$ by
\[   D_{\beta} \psi  = \sum_{i=1}^n  \beta^{-1} (E_i) \cdot  \nabla_{E_i} \psi  = \sum_{i=1}^n  E_i \cdot  \nabla_{\beta^{-1}(E_i)}  \psi  ,  \]
where $\beta$ was identified with the induced (1,1)-tensor $\beta$
via $\beta(X, Y) = g(X, \beta(Y))$.
Recall that a symmetric (0,
2)-tensor field $\beta$ is called  \\
(i) a {\it divergencefree tensor} if $ {\rm div}(\beta) =
\sum_{i=1}^n ( \nabla_{E_i} \beta)(E_i) = 0$.     \\
(ii)  a {\it Codazzi tensor} if $(\nabla_X \beta)(Y, Z) = (\nabla_Y
\beta)(X, Z)$ holds for all vector fields $X, Y, Z$.

\bigskip
Let $( \ , \  ) := {\rm Re} \langle \ , \  \rangle$ denote the real
part of the standard Hermitian  product $\langle \ , \ \rangle$ on
the spinor bundle $\Sigma (M)$ over $M^n$. Let $\alpha$ be a 1-form
on $M^n$ induced by a nondegenerate symmetric tensor $\beta$ and
spinor fields $\phi, \psi  \in \Gamma(\Sigma)$ via
\[  \alpha(X) = ( \phi , \ \beta^{-1}(X) \cdot \psi ) .   \]
Then
\[  {\rm div}(\alpha) = - ( D_{\beta} \phi, \, \psi ) + ( \phi, \,
D_{\beta} \psi ) + ( \phi,  \, {\rm div}( \beta^{-1}) \cdot \psi ) .
\]  Consequently, if $\beta^{-1}$ is divergencefree,  then  $D_{\beta}$ is a self-adjoint
elliptic operator of first order and hence its spectrum is discrete and real.

\bigskip
We have proved in [6, 7] the following theorems.

\begin{thm}
 Let $(M^n, g)$ be an n-dimensional closed
Riemannian spin manifold.  Let $\beta$ be such a nondegenerate
symmetric tensor on $M^n$  that both ${\rm div} ( \beta^{-1} ) = 0$
and ${\rm tr}(\beta^{-1}) = 0$ vanish identically.    Let $\lambda_1
 \in {\mathbb R} $ and $\overline{\lambda}_1 \in {\mathbb R} $
be the first eigenvalue of $D$ and $D_{\beta}$, respectively. Then
we have
\[
\lambda_1^2  \,   \geq   \,   \inf_M  \Bigg\{  \frac{n \, S}{4(n-1)}
+ \frac{n \, \overline{\lambda}_1^2  }{(n-1) \, \vert \beta^{-1}
\vert^2 } + \frac{ n \, \triangle ( \vert \beta^{-1}  \vert^2 )}{2
(n-1)\vert \beta^{-1}  \vert^2}   \Bigg\}.
\]

\bigskip
\indent The limiting case  occurs  if and only if  there exists a
spinor field $\psi_1$ on $(M^n, g)$ with the following
properties:   \\
$(i)$ The differential equation
\[  \nabla_X \psi_1 =   - \frac{\lambda}{n}  \, X
\cdot \psi_1  - \frac{\overline{\lambda}}{\vert \beta^{-1} \vert^2}
\, \beta^{-1}(X) \cdot \psi_1   \]
 holds for some
constants $\lambda , \, \overline{\lambda} \in {\mathbb R}$ and for
all vector fields $X$.   \\
$(ii)$ $\psi_1$ is a first eigenspinor of both $D$ and $D_{\beta}$.
\end{thm}

\bigskip
\begin{thm}
Let $(M^n, g)$ be an n-dimensional
closed Riemannian spin manifold and consider a nondegenerate Codazzi
tensor $\beta$ such that ${\rm tr}(\beta^{-1}) = 0$ vanishes
identically. Denote by $\overline{g}$ the metric induced by $\beta$
via $\overline{g} (X, Y) = g ( \beta(X), \beta(Y) )$ and by
$\overline{D}$ the Dirac operator of $\overline{g}$.
 Let $\lambda_1  \in {\mathbb R} $
and $\overline{\lambda}_1 \in {\mathbb R} $ be the first eigenvalue
of the Dirac operators $D$ and $\overline{D}$, respectively.
 Then we have
\[
\lambda_1^2  \,   \geq   \,   \inf_M  \Bigg\{  \frac{n \, S}{4(n-1)}
+ \frac{n \, \overline{\lambda}_1^2  }{(n-1) \, \vert \beta^{-1}
\vert^2 } + \frac{ n \, \triangle ( \vert  {\rm det} ( \beta^{-1} )
\vert  \,  \vert \beta^{-1} \vert^2 ) }{2 (n-1) \vert  {\rm det} (
\beta^{-1} ) \vert  \,  \vert \beta^{-1} \vert^2}   \Bigg\}.
\]

\bigskip
\indent The limiting case of occurs  if and only if there exists a
spinor field $\psi_1$ on $(M^n, g)$ with the following
properties:   \\
$(i)$ The differential equation
 \[ \nabla_X \psi_1 =   -
\frac{\lambda}{n} \, X \cdot \psi_1 -
\frac{\overline{\lambda}}{\vert \beta^{-1} \vert^2} \, \beta^{-1}(X)
\cdot \psi_1   \] holds for some constants $\lambda , \,
\overline{\lambda}  \in {\mathbb R}$ and for all vector fields $X$.
\\
$(ii)$ $\psi_1$ is a first eigenspinor of both $D$ and
$\overline{D}$.
\end{thm}

\noindent
\section{Dirac eigenvalues estimates over Sasakian manifolds}

\noindent In this section we will apply Theorem 2.1 to Sasakian manifolds.  Consider a manifold $M^{2m+1}$  of odd dimension $n=2m+1$. An
almost contact metric structure $( \phi, \xi, \eta, g)$ of
$M^{2m+1}$ consists of a (1,1)-tensor field $\phi$, a vector field
$\xi$, a 1-form $\eta$, and a metirc $g$ with the following
properties:
\[   \eta (\xi) = 1, \qquad  \phi^2(X) = -X + \eta(X) \xi,  \qquad g (\phi X, \, \phi Y ) = g(X, Y) - \eta(X) \eta(Y) . \]
The {\it fundamental $2$-form} $\Phi$ of the contact structure is a 2-form defined by
\[ \Phi ( X, Y) = g ( X, \, \phi(Y) )  .  \]
An almost contact metric structure $( \phi, \xi, \eta, g)$ of
$M^{2m+1}$ becomes a Sasakian structure if
\[   ( \nabla_X \phi )(Y) = g(X, Y) \xi - \eta(Y) X    \]
holds for all vector fields $X, Y$.
A Sasakian manifold $( M^{2m+1}, \phi, \xi, \eta, g)$ is called {\it
eta-Einstein}  if the Ricci curvature tensor ${\rm Ric}$
satisfies
\begin{equation}
{\rm Ric} = \kappa \, g +  \tau  \eta \otimes \eta
\end{equation}
for some constants $\kappa , \tau  \in {\mathbb R}$ with $\kappa +
\tau = 2m$.  Any eta-Einstein Sasakian manifold is necessarily of
constant scalar curvature $S$ and we can rewrite eta-Einstein
condition (3.1) as
\[
{\rm Ric} = \left( \frac{S}{n-1} - 1 \right)  g +  \left( n -
\frac{S}{n-1} \right)  \eta \otimes \eta ,  \qquad  n = 2m +1 .
\]

\bigskip
From now on we assume that any Sasakian manifold $( M^{2m+1}, \phi, \xi, \eta, g)$ we consider has a fixed spin structure.
An important property of a Sasakian spin manifold $(
M^{2m+1}, \phi, \xi, \eta, g)$ is that the spinor bundle $\Sigma(M)$  splits under the action of the
fundamental 2-form $\Phi$ as follows.

\begin{lem}
Let $( M^{2m+1}, \phi, \xi, \eta, g)$ be an almost contact metric
manifold with spin structure and fundamental $2$-form $\Phi$.  Then
the spinor bundle
 $\Sigma$ splits into the orthogonal direct sum $\Sigma = \Sigma_0 \oplus \Sigma_1 \oplus \cdots \oplus \Sigma_m$ with    \\
(i) $\Phi \vert_{\Sigma_r} =  \sqrt{-1} (2r - m) I,  \qquad {\rm
dim}(\Sigma_r) =  {m \choose r} \qquad (0 \leq r \leq m),
 $   \\
(ii) $ \xi \vert_{\Sigma_0 \oplus \Sigma_2 \oplus \Sigma_4 \oplus
\cdots}  = ( \sqrt{-1} )^{2m+1} I,  \qquad
\xi \vert_{\Sigma_1 \oplus \Sigma_3 \oplus \Sigma_5 \oplus \cdots}  = - ( \sqrt{-1} )^{2m+1} I, $    \\
where $I$ stands for the identity map.
 Moreover, the bundles $\Sigma_0$ and $\Sigma_m$ can be defined
by
\begin{eqnarray*}
\Sigma_0 & = & \{ \,  \psi \in \Sigma \ : \ \phi(X) \cdot \psi + \sqrt{-1} X \cdot \psi + (-1)^m \eta(X) \psi = 0 \ \mbox{for all vectors} \, X  \, \} ,   \\
\Sigma_m & = & \{ \,  \psi \in \Sigma \ : \ \phi(X) \cdot \psi - \sqrt{-1} X \cdot \psi - \eta(X) \psi = 0 \ \mbox{for all vectors} \, X  \, \}.
\end{eqnarray*}
In particular, we have the formulas
\begin{eqnarray*}
&     &  \xi \cdot \psi_0 = (-1)^m \sqrt{-1} \psi_0 ,  \qquad  \Phi \cdot \psi_0 = - m \sqrt{-1} \psi_0 , \qquad \psi_0 \in  \Sigma_0 ,       \\
&     &  \xi \cdot \psi_m =  \sqrt{-1} \psi_m ,  \qquad  \Phi \cdot \psi_m =  m \sqrt{-1} \psi_m , \qquad \psi_m \in  \Sigma_m .
\end{eqnarray*}
\end{lem}

\bigskip
Over Sasakian spin  manifolds, a special class of spinors
deserves attention.

\bigskip \noindent {\bf Definition 3.1}
  A nontrivial spinor field
$\psi$ on Sasakian spin manifold $( M^{2m+1}, \phi, \xi, \eta, g)$
is called an {\it eta-Killing spinor} with Killing pair $(a, b)$ if
it satisfies
\begin{equation}  \nabla_X \psi = a \, X \cdot \psi + b \eta(X) \xi  \cdot \psi
\end{equation}
for some real numbers $a , b \in {\mathbb R}, \ a \not= 0,$ and for
all vector fields $X$.

\bigskip
\noindent Note that if $b =0$, then equation (3.2) reduces to equation (1.2). Moreover, any eta-Killing spinor with Killing pair $(a, b)$ is an
eigenspinor of the Dirac operator with eigenvalue $\lambda = -
(2m+1)a - b$.

\bigskip
Now we summarize some basic relations between
the Killing pair $(a, b)$ of an eta-Killing spinor and the geometry
of the Sasakian manifold. For proofs for Propositions 3.1-3.4 we refer to [5, 7].
In the following we will often write $n$ to mean the
dimension $2m+1$ of the manifold $M^{2m+1}$.

\begin{pro}
Let $( M^{2m+1}, \phi, \xi, \eta, g) ,  \,  m  \geq 2  ,$ be a
Sasakian spin manifold and suppose that it admits an eta-Killing
spinor $\psi$ with Killing pair $(a, b)$, where both $a \not=0$ and $b \not=
0$ are nonzero. Then $( M^{2m+1}, \phi, \xi, \eta, g)$ is
eta-Einstein with scalar curvature $S= 4n(n-1) a^2 + 8(n-1) ab$. Moreover, all the possible values for $a, b$ can be
expressed in terms of the scalar curvature as
\[  ( a, b) \ = \  \left( \frac{1}{2} , \  - \frac{n}{4} +
\frac{S}{4(n-1)} \right),  \qquad     \left( - \frac{1}{2} , \
\frac{n}{4} - \frac{S}{4(n-1)} \right) ,    \]
and the following statements are true: \\
(i) If $(a, b) = \left( \frac{1}{2} , \  - \frac{n}{4} +
\frac{S}{4(n-1)} \right)$, then $m \equiv 0 \ {\rm mod} \, 2$ and
$\psi \in \Gamma ( \Sigma_0 )$ is a section in $\Sigma_0$.  \\
(ii) If $(a, b) = \left( - \frac{1}{2} , \   \frac{n}{4} -
\frac{S}{4(n-1)} \right)$ and $m \equiv 0 \ {\rm mod} \, 2$, then
$\psi \in \Gamma ( \Sigma_m )$ is a section in $\Sigma_m$.   \\
(iii) If $(a, b) =  \left( - \frac{1}{2} , \   \frac{n}{4} -
\frac{S}{4(n-1)} \right)$ and $m \equiv 1 \ {\rm mod} \, 2$, then
$\psi \in  \Gamma ( \Sigma_0 ) \cup \Gamma ( \Sigma_m )$ is a
section in $\Sigma_0$ or in $\Sigma_m$.
\end{pro}

\begin{pro}
Let $( M^{2m+1}, \phi, \xi, \eta, g) ,  \,  m  \geq 2  ,$ be a
simply-connected Sasakian spin manifold. Suppose that $( M^{2m+1},
\phi, \xi, \eta, g)$ is eta-Einstein.  Then,  in case   \\
(i) \,  $m \equiv 0 \ {\rm mod} \, 2$,  there exists an
eta-Killing spinor $\psi_0 \in  \Gamma ( \Sigma_0 )$ with Killing pair $(
\frac{1}{2} , \,  - \frac{n}{4} + \frac{S}{4(n-1)} )$ as well as an
eta-Killing spinor $\psi_m \in  \Gamma ( \Sigma_m )$ with Killing pair $( -
\frac{1}{2} , \,   \frac{n}{4} - \frac{S}{4(n-1)} )$.    \\
(ii) \,  $m \equiv 1 \ {\rm mod} \, 2$, there exist two
eta-Killing spinors $\psi_0, \, \psi_m$  with Killing pair $( - \frac{1}{2} ,
\,
 \frac{n}{4} - \frac{S}{4(n-1)} )$  such that $\psi_{\alpha}$ is a
 section in the bundle $\Sigma_{\alpha} \ ( \alpha = 0 , m ).$
\end{pro}

\begin{pro}
Let $( M^3, \phi, \xi, \eta, g) $ be a  $3$-dimensional Sasakian
spin manifold and suppose that it admits an eta-Killing spinor
$\psi$ with Killing pair $(a, b)$, where $a \not=0$ and $b \not=0$.
Then $( M^3, \phi, \xi, \eta, g)$ is eta-Einstein with constant
scalar curvature $S= 24 a^2 + 16ab$. Moreover, all the possible
values for $a, b$ can be expressed in terms of the scalar curvature
as
  \begin{eqnarray*} (a,
b) & = & \left( - \frac{1}{2} , \ \frac{3}{4} - \frac{S}{8} \right)
, \ \left( \frac{-2 + \sqrt{4 + 2S}}{4} , \  \frac{4 -
\sqrt{4 + 2 S}}{4} \right) , \\
&     &  \left( \frac{-2 - \sqrt{4 + 2S}}{4} , \  \frac{4 + \sqrt{4
+ 2 S}}{4}  \right).
\end{eqnarray*}
\end{pro}

\begin{pro}
Let $( M^3, \phi, \xi, \eta, g)$ be a simply-connected Sasakian spin
manifold of dimension $3$ and suppose that the scalar curvature $S$
of $g$
is constant. Then,   \\
(i) \,  there exist two eta-Killing spinors $\psi_0, \, \psi_1$ with Killing pair $( - \frac{1}{2} , \,
 \frac{3}{4} - \frac{S}{8} )$  such that $\psi_{\alpha}$ is a
 section in the bundle $\Sigma_{\alpha} \ ( \alpha = 0 , 1 ).$    \\
 (ii) \, If $S \geq - 2$, there exists an eta-Killing spinor $\psi  \in
\Gamma( \Sigma = \Sigma_0 \oplus \Sigma_1 )$ with Killing pair
$\left( \frac{-2 \pm \sqrt{4 + 2S}}{4} , \  \frac{4 \mp \sqrt{4 + 2
S}}{4}  \right)$.
\end{pro}

\bigskip
We are now ready to apply Theorems 2.1 to Sasakian
spin manifolds. The resulting inequality (3.3) clearly improves inequality (1.1).

\begin{pro}
Let $( M^{2m+1}, \phi, \xi, \eta, g) ,  \,  m  \geq 1  ,$ be a
closed Sasakian spin manifold. Let $\beta^{-1}$ be a nondegenerate
symmetric tensor field on $M^{2m+1}$ defined by $\beta^{-1} =
\frac{2}{n} \, I - 2 \, \xi \otimes \eta$. (Note that ${\rm div}(
\beta^{-1} ) = 0$ and ${\rm tr} ( \beta^{-1} ) = 0$ .) Let
$\lambda_1   \in {\mathbb R} $ and $\overline{\lambda}_1 \in
{\mathbb R} $ be the first eigenvalue of $D$ and $D_{\beta}$,
respectively. Then we have
\begin{equation}
\lambda_1^2  \,   \geq   \,      \frac{n \, S_{\rm min}}{4(n-1)}
 + \frac{n^2 \, \overline{\lambda}_1^2  }{4(n-1)^2}  ,
\end{equation}
where $ S_{\rm min}$ denotes the minimum of the scalar curvature.
 The limiting case of $(3.3)$ occurs, in case    \\
(i)  $n \geq 5$, if and only if there exists an eta-Killing spinor
$\psi_1$ with Killing pair
\begin{equation}
\left( \frac{1}{2} , \  - \frac{n}{4} +
\frac{S}{4(n-1)} \right),  \qquad    \left( - \frac{1}{2} , \
\frac{n}{4} - \frac{S}{4(n-1)} \right) ,
\end{equation}
such that $\psi_1$ is a
first eigenspinor of both $D$ and $D_{\beta}$.  \\
 (ii) $n =3$, if and only if
there exists an eta-Killing spinor $\varphi_1$ with Killing pair
\begin{equation}
  \left( \frac{-2 + \sqrt{4 + 2S}}{4} , \  \frac{4 - \sqrt{4
+ 2 S}}{4}  \right)
\end{equation}
 such that $\varphi_1$ is a
first eigenspinor of both $D$ and $D_{\beta}$.  \\
\end{pro}

\bigskip
Let $( M^{2m+1}, \phi, \xi, \eta, g), \, m \geq 1,$
be a closed Sasakian spin manifold with positive scalar curvature $S > 0$. From inequality (3.3) we see
that the first eigenvalue $\lambda_1 \not=0$ is necessarily nonzero. The statement for the limiting case of (3.3) then
  gives rise to a natural question:

\bigskip \noindent
 Is every eta-Killing spinor with Killing pair (3.4) or (3.5) a first eigenspinor of the Dirac operator
 ?

\bigskip \noindent
We have recently found that  answer to the question in 3-dimensional case is positive [8], but the question in higher dimensional case is still open.

\begin{Literature}{xx}
\bibitem{1}
C. B\"{a}r,  Real Killing spinors and holonomy,  Commun. Math.
Phys. 154 (1993) 509-521.
\bibitem{2}
H. Baum, Th. Friedrich, R. Grunewald, I. Kath,  Twistors and
Killing spinors on Riemannian manifolds, Teubner,
Leipzig/Stuttgart, 1991.
\bibitem{3}
Th. Friedrich,  Der erste Eigenwert des Dirac-Operators einer
kompakten Riemannschen Mannigfaltigkeit nichtnegativer
Skalarkr\"{u}mmung,  Math. Nachr. 97  (1980)  117-146.
\bibitem{4}
Th. Friedrich, Dirac operators in Riemannian geometry, Graduate studies in mathematics Vol. 25, AMS, Providence, 2000.
\bibitem{5}
Th. Friedrich, E.C. Kim,  The Einstein-Dirac equation on
Riemannian spin manifolds,  J. Geom. Phys. 33 (2000) 128-172.
\bibitem{6}
Th. Friedrich, E.C. Kim, Eigenvalues estimates for the Dirac
operator in terms of Codazzi tensors, Bull. Korean Math. Soc. 45
(2008) 365-373.
\bibitem{7}
E.C. Kim, Dirac eigenvalues estimates  in terms of divergencefree
symmetric tensors, to appear in Bull. Korean Math. Soc.
\bibitem{8}
E.C. Kim, The first eigenvalue of the Dirac operator on non-Einstein
Sasakian 3-manifolds, in preparation.
\bibitem{9}
K.-D. Kirchberg, An estimation for the first eigenvalue of the
Dirac operator on closed K\"{a}hler manifolds of positive scalar
curvature, Ann. Global Anal. Geom. 4 (1986) 291-325.
\bibitem{10}
K.-D. Kirchberg, U. Semmelmann, Complex Contact Structures and the
First Eigenvalue of the Dirac Operator on Kahler Manifolds, Geom.
and Funct. Analysis 5 (1995) 604-618.
\bibitem{11}
W. Kramer, U. Semmelmann, G. Weingart, {\it Eigenvalue estimates for
the Dirac operator on quaternionic K\"{a}hler manifolds}, Math. Z.
230 (1999) 727-751.
\end{Literature}

\bigskip   \noindent
Department of Mathematics Education, \\   Andong National
University,
Andong 760-749, South  Korea     \\
e-mail: eckim@andong.ac.kr

\end{document}